\documentclass[12pt]{article}

\usepackage[cp850]{inputenc}
\usepackage{amssymb}
\usepackage{latexsym}

\topmargin=-\headheight
\def\Z{\mathbb{Z}}
\def\R{\mathbb{R}}
\def\C{\mathbb{C}}
\def\S{{\cal S}}
\def\A{{\bf A}}
\def\B{{\bf B}}
\def\b{{\bf b}}

\def\t{{\bf t}}
\def\x{{\bf x}}

\begin{document}
\setlength{\parindent}{0pt}
\setlength{\parskip}{0.4cm}

\newtheorem{theorem}{Theorem}
\newtheorem{lemma}[theorem]{Lemma}
\newtheorem{definition}{Definition}

\begin{center}

\large{\bf A Closer Look at Lattice Points in Rational Simplices}
\footnote{This work is part of the author's Ph.D. thesis. \\ 
           {\it Mathematical Reviews Subject Numbers:} 05A15, 11D75.}

\vspace{5mm}

\normalsize

{\sc Matthias Beck} \footnote{\tt 
http://www.math.temple.edu/$\sim$matthias} \\
Dept. of Mathematics, Temple University \\
Philadelphia, PA 19122 \\
{\tt matthias@math.temple.edu}

Submitted: March 19, 1999; accepted: September 14, 1999

\end{center}

\vspace{5mm}

\footnotesize
{\bf Abstract.} We generalize Ehrhart's idea (\cite{ehrhart}) of counting 
lattice points in
dilated rational polytopes: Given a rational simplex, that is, an 
$n$-dimensional polytope
with $n+1$ rational vertices, we use its description as the intersection of 
$n+1$ halfspaces,
which determine the facets of the simplex. Instead of just a single 
dilation factor, we allow
different dilation factors for each of these facets. We give an elementary 
proof that the
lattice point counts in the interior and closure of such a {\it 
vector-dilated} simplex
are quasipolynomials satisfying an Ehrhart-type reciprocity law. This 
generalizes the classical
reciprocity law for rational polytopes (\cite{macdonald}, \cite{mcmullen}, 
\cite{stanley}).
As an example, we derive a lattice point count formula for a rectangular 
rational
triangle, which enables us to compute the number of lattice points inside 
any rational polygon.
\normalsize


\section{\normalsize Introduction}
One of the exercises on the greatest integer function $ [ x ] $ in an 
elementary course in
Number Theory is to prove the statement
   \begin{equation}\label{gi} \left[ \frac{ t-1 }{ a }  \right] = - \left[ 
\frac{ -t }{ a }  \right] - 1 \end{equation}
for any integers $ t, a \not= 0$.
Geometrically, this is a special instance of a much more general theme. 
Consider the interval
$ \left[ 0 , \frac{ 1 }{ a }  \right] $, viewed as a 1-dimensional rational
polytope. (A {\it rational polytope} is a polytope whose vertices are 
rational.) Now we dilate this polytope by
an integer factor $t>0$, and count the number of integer points (''lattice 
points'') in the
dilated polytope. It is straightforward that this number in the open 
dilated polytope is
$ \left[ \frac{ t-1 }{ a }  \right] $, whereas in the closure there are
$ \left[ \frac{ t }{ a }  \right] + 1 $ integer points.

More generally, let ${\cal P}$ be an $n$-dimensional convex rational 
polytope in $ \R^{ n } $.
For $ t \in \Z_{ >0 }  $, let $ L ( {\cal P}^{\circ}  , t ) = \# \left( t 
{\cal P}^{\circ}  \cap \Z^{ n }  \right) $
and $ L ( \overline{\cal P} , t ) = \# \left( t \overline{\cal P} \cap \Z^{ 
n }  \right) $ be
the number of lattice points in the interior of the dilated polytope $ t 
{\cal P} = \{ tx : x \in {\cal P}  \} $
and its closure, respectively. That is, if ${\cal P}$ denotes the above 
1-dimensional polytope, we have
   \[ L ( {\cal P}^{\circ}  , t ) = \left[ \frac{ t-1 }{ a }  \right] \quad 
\mbox{ and } \quad L ( \overline{\cal P} , t ) = \left[ \frac{ t }{ a 
}  \right] + 1 \ . \]
There are two remarkable features hidden in these expressions: First, we have
   \begin{theorem}\label{quasi} $ L ( {\cal P}^{\circ}  , t ) $ and $ L ( 
\overline{\cal P} , t ) $ are quasipolynomials in $t$. \end{theorem}
A {\it quasipolynomial} is an expression of the form $ c_{n}(t) \ t^{n} + 
\dots + c_{1}(t) \ t + c_{0}(t) $,
where $c_{0}, \dots , c_{n}$ are periodic functions in $t$. Theorem 
\ref{quasi} is easily
verified for our one-dimensional polytope by writing
$ [x] = x - \{ x \} $, where $ \{ x \} $ denotes the fractional part of $x$.
Moreover, viewing both these quasipolynomials as algebraic expressions in 
the {\it integer}
variable $t$, (\ref{gi}) becomes a reciprocity law:
   \begin{theorem}\label{recipr} $ L ( {\cal P}^{\circ}  , -t ) = (-1)^{ n 
} L ( \overline{\cal P} , t ) $. \end{theorem}
Both Theorem \ref{quasi} and \ref{recipr} are true for any rational 
polytope ${\cal P}$.
The proof of Theorem \ref{quasi} is due to
Ehrhart, who initiated the study of the lattice point count in dilated 
polytopes (\cite{ehrhart}).
He conjectured Theorem \ref{recipr}, which was first proved by Macdonald
(for the case that ${\cal P}$ has integer vertices, \cite{macdonald}), 
later also by McMullen
(\cite{mcmullen}), and Stanley (\cite{stanley}).

We generalize the notion of dilated polytopes for rational simplices, that 
is, rational polytopes of dimension
$n$ with $n+1$ vertices. We use the description of a simplex as the 
intersection of $n+1$ halfspaces,
which determine the facets of the simplex: Instead of dilating the simplex 
by a single factor,
we allow different dilation factors for each facet.
\begin{definition}\label{def} Let the rational simplex $ \S_{ \A } $ be 
given by
   \[ \S_{ \A } = \left\{ \x \in \R^{ n } : \ \A \ \x \leq \b \right\} \ , \]
with $ \A \in M_{ (n+1) \times n } ( \Z ) , \b \in \Z^{ n+1 } $. Here the 
inequality is
understood component\-wise. For $ \t \in \Z^{ n+1 } $, define the {\bf 
vector-dilated simplex}
$ \S_{ \A }^{ ( \t ) }  $ as
   \[ \S_{ \A }^{ ( \t ) } = \left\{ \x \in \R^{ n } : \ \A \ \x \leq \t 
\right\} \ . \]
For those $ \t $ for which $ { S_{ \A }^{ ( \t ) }   } $ is nonempty and 
bounded,
we define the number of lattice points in the interior and closure of
$ \S_{ \A }^{ ( \t ) }  $ as
   \[ L \left( \S_{ \A }^{\circ}  , \t \right) = \# \left( \left( \S_{ \A 
}^{ ( \t ) }  \right)^{\circ}  \cap \Z^{ n }  \right) \quad \mbox{ and } 
\quad L \left( \overline{\S_{ \A }} , \t \right) = \# \left( {\S_{ \A }^{ ( 
t ) }} \cap \Z^{ n }  \right) \ , \]
respectively.
\end{definition}
Geometrically, we fix for a given simplex the normal vectors to its facets 
and consider all
possible positions of these normal vectors that 'make sense'.
The previously defined quantities $ L ( {\cal P}^{\circ}  , t ) $ and $ L ( 
\overline{\cal P} , t ) $
can be recovered from this new definition by choosing $ \t = t \b $.
The corresponding result to Theorems \ref{quasi} and \ref{recipr} is
\begin{theorem}\label{multithm}
$ L \left( \S_{ \A }^{\circ}  , \t \right) $ and $ L \left( \overline{\S_{ 
\A }} , \t \right) $ are quasipolynomials
in $ \t \in \Z^{ n+1 } $, satisfying
   \begin{equation}\label{rec} L \left( \S_{ \A }^{\circ}  , - \t \right) = 
(-1)^{ n }  L \left( \overline{\S_{ \A }} , \t \right) \ . \end{equation}
\end{theorem}
A {\it quasipolynomial} in the $d$-dimensional variable $\t$ is the obvious 
generalization of
a quasipolynomial in a 1-dimensional variable.

We give an elementary proof of Theorem \ref{multithm}, only relying on 
(\ref{gi}) and a basic
lemma on quasipolynomials. Theorems \ref{quasi} and \ref{recipr} follow as 
immediate
corollaries, considering the fact that any polytope can be triangulated 
into simplices.
In fact, the original motivation for Theorem \ref{multithm} was to 
construct an elementary
proof of Theorem \ref{recipr}.


\section{\normalsize A lemma on quasipolynomials}
\begin{lemma}\label{quasilemma} Let $ q( t_{ 1 }, \dots , t_{ m }  ) $ be a 
quasipolynomial, and
fix $ a_{ 1 } , \dots , a_{ m } $, $ c_{ 0 } , \dots , c_{ m } $, $ d \in 
\Z , d \not= 0 $. Then
   \[ Q_{ 1 } ( \t ) = Q_{ 1 } ( t_{ 0 }, t_{ 1 } , \dots , t_{ m }  ) = 
\sum_{ k=1 }^{ \left[ \frac{ c_{ 0 } t_{ 0 } + \dots + c_{ m } t_{ m }  - 1 
}{ d }  \right]  } q \left( t_{ 1 } + a_{ 1 } k , \dots , t_{ m } + a_{ m } 
k \right)   \]
and
   \[ Q_{ 2 } ( \t ) = \sum_{ k=0 }^{ \left[ \frac{ c_{ 0 } t_{ 0 } + \dots 
+ c_{ m } t_{ m } }{ d } \right]  } q \left( t_{ 1 } + a_{ 1 } k , \dots , 
t_{ m } + a_{ m } k \right)   \]
are also quasipolynomials.
\end{lemma}
{\it Remark.} Here and in the following we define a finite series $ \sum_{ 
k=a }^{ b } \dots $
for {\it both} cases $ a \leq b $ and $ a > b $, in the usual way:
   \begin{equation}\label{sum} \sum_{ k=a }^{ b } \dots = \left\{ 
\begin{array}{ll} \sum_{ k=a }^{ b } \dots & \mbox{ if } a \leq b \\
                                                           0 & \mbox{ if } 
a = b+1 \\
                                                           - \sum_{ k=b+1 
}^{ a-1 } \dots & \mbox{ if } a \geq b+2 \end{array}\right.  \end{equation}
{\it Proof.}
We will prove the statement for $ Q_{ 2 } $; the proof for $ Q_{ 1 } $ 
follows in a similar fashion.
After writing $q$ in all its terms and multiplying out the binomial 
expressions, it
suffices to prove that
   \[ Q_{ 3 } ( \t ) = \sum_{ k=0 }^{ \left[ \frac{ c_{ 0 } t_{ 0 } + \dots 
+ c_{ m } t_{ m } }{ d } \right] } f \left( t_{ 1 } + a_{ 1 } k , \dots , 
t_{ m } + a_{ m } k \right) \ k^{ j }   \]
is a quasipolynomial, where $j$ is a fixed nonnegative integer and $f$ is a 
periodic function
in $m$ variables. Consider a period $p$ which is common to all the 
arguments of $f$, that is,
$ f \left( x_{ 1 } + p , \dots , x_{ m } + p \right)  = f \left( x_{ 1 } , 
\dots , x_{ m }  \right)  $.
To see that $ Q_{ 3 } $ is a quasipolynomial, use the properties of $f$ to 
write it as
\begin{eqnarray*} &\mbox{}& Q_{ 3 } ( \t ) = f \left( t_{ 1 } , \dots , t_{ 
m } \right) \ \sum_{ k=0 }^{ \left[ \frac{ c_{ 0 } t_{ 0 } + \dots + c_{ m 
} t_{ m } }{ d p } \right] } (kp)^{ j } \\
                   &\mbox{}& \qquad + \ f \left( t_{ 1 } + a_{ 1 } , \dots 
, t_{ m } + a_{ m } \right) \ \sum_{ k=0 }^{ \left[ \frac{ c_{ 0 } t_{ 0 } 
+ \dots + c_{ m } t_{ m } - d }{ d p } \right] } (1+kp)^{ j } \  + \\
                   &\mbox{}& \qquad + \ f \left( t_{ 1 } + 2 a_{ 1 } , 
\dots , t_{ m } + 2 a_{ m } \right) \ \sum_{ k=0 }^{ \left[ \frac{ c_{ 0 } 
t_{ 0 } + \dots + c_{ m } t_{ m } - 2d }{ d p } \right] } (2+kp)^{ j } \ + 
\ \dots \ + \\
                   &\mbox{}& \qquad + \ f \Bigl( t_{ 1 } + (p-1) a_{ 1 } , 
\dots , t_{ m } + (p-1) a_{ m } \Bigr) \ \sum_{ k=0 }^{ \left[ \frac{ c_{ 0 
} t_{ 0 } + \dots + c_{ m } t_{ m } - (p-1) d }{ d p } \right] } (p-1+kp)^{ 
j } \ . \end{eqnarray*}
Upon expanding all the binomials, putting the finite sums into closed forms,
and writing $[x] = x - \{ x \} $, the only dependency on $\t$ is periodic 
(with period
dividing $ d p $) or polynomial.
\hfill {} $\Box$


\section{\normalsize Proof of Theorem \ref{multithm} }
We induct on the dimension $n$. First, a 1-dimensional rational simplex 
${\S_{ \A }}$ is an interval with
rational endpoints. Hence $ \S_{ \A }^{ (\t) } $ is given by
   \[ \frac{ t_{ 1 }  }{ a_{ 1 }  } \leq x \leq \frac{ t_{ 2 }  }{ a_{ 2 
}  } \ ,  \]
so that we obtain
   \[ L \left( \S_{ \A }^{\circ}  , \t \right) = \left[ \frac{ t_{ 2 } - 1 
}{ a_{ 2 }  }  \right] - \left[ \frac{ t_{ 1 } }{ a_{ 1 
}  }  \right]  \quad \mbox{ and } \quad L \left( \overline{\S_{ \A }} , \t 
\right) = \left[ \frac{ t_{ 2 } }{ a_{ 2 }  }  \right] - \left[ \frac{ t_{ 
1 } - 1 }{ a_{ 1 }  }  \right]   \ . \]
These are quasipolynomials, as can be seen, again, by writing $[x] = x - \{ 
x \} $. Furthermore, by (\ref{gi}),
   \[ L \left( \S_{ \A }^{\circ}  , - \t \right) = \left[ \frac{ - t_{ 2 } 
- 1 }{ a_{ 2 }  }  \right] - \left[ \frac{ - t_{ 1 } }{ a_{ 1 }  }  \right] 
= - \left[ \frac{ t_{ 2 } }{ a_{ 2 }  }  \right] + \left[ \frac{ t_{ 1 } - 
1 }{ a_{ 1 }  }  \right]  = - L \left( \overline{\S_{ \A }} , \t \right) \ . \]
Now, let $\S_{ \A }$ be an $n$-dimensional rational simplex. After harmless 
unimodular
transformations, which leave the lattice point count invariant, we may 
assume that the defining
inequalities for $\S_{ \A }$ are
\begin{eqnarray*}\begin{array}{rllllll} a_{ 11 } x_{ 1 } & & & & & \leq & 
b_{1} \\
                   a_{ 21 } x_{ 1 } &+& \dots &+& a_{ 2n } x_{ n } & \leq & 
b_{2} \\
                   & \vdots \\
                   a_{ n+1,1 } x_{ 1 } &+& \dots &+& a_{ n+1,n } x_{ n } & 
\leq & b_{n+1} \ . \end{array} \end{eqnarray*}
(Actually, we could obtain an lower triangular form for $ \A $; however, 
the above form suffices for our purposes.)
Hence there exists a vertex $ {\bf v} = ( v_{ 1 } , \dots , v_{ n } ) $ 
with $ v_{ 1 } = \frac{ b_{ 1 }  }{ a_{ 11 }  } $
and another vertex $ {\bf w} = ( w_{ 1 } , \dots , w_{ n } ) $ whose first 
component is not $ \frac{ b_{ 1 }  }{ a_{ 11 }  } $.
After switching $ x_{ 1 } $ to $ -x_{ 1 } $, if necessary,
we may further assume that $ v_{ 1 } < w_{ 1 } $. Since $ {\bf w} $ 
satisfies all
equalities but the first one, it is not hard to see that $ {\bf w} $ has 
first component
$ w_{1} = r_{ 2 } b_{ 2 } + \dots + r_{ n } b_{ n } $ for some rational 
numbers $ r_{ 2 } , \dots , r_{ n } $;
write this number as $ w_{1} = \frac{ c_{ 2 } b_{ 2 } + \dots + c_{ n } b_{ 
n }  }{ d } $
with $ c_{ 2 } , \dots , c_{ n } , d \in \Z $.
Viewing the defining inequalities of the vector-dilated simplex $ \S_{ \A 
}^{ ( \t ) } $ as
\begin{eqnarray*}\begin{array}{rllllll} & & \frac{ t_{ 1 } }{ a_{ 11 }  } & 
\leq & x_{ 1 } & \leq & \frac{ c_{ 2 } t_{ 2 } + \dots + c_{ n } t_{ n 
}  }{ d }  \\
                   a_{ 22 } x_{ 2 } &+& \dots &+& a_{ 2n } x_{ n } & \leq & 
t_{ 2 } - a_{ 21 } x_{ 1 }  \\
                   & \vdots \\
                   a_{ n+1,2 } x_{ 2 } &+& \dots &+& a_{ n+1,n } x_{ n } & 
\leq & t_{ n+1 } - a_{ n+1,1 } x_{ 1 }  \ , \end{array} \end{eqnarray*}
we can compute the number of lattice points in the interior and closure of 
$ \S_{ \A }^{ ( \t ) } $ as
   \begin{equation}\label{int} L \left( \S_{ \A }^{\circ}  , \t \right) = 
\sum_{ m = \left[ \frac{ t_{ 1 } }{ a_{ 11 }  }  \right] + 1 }^{ \left[ 
\frac{ c_{ 2 } t_{ 2 } + \dots + c_{ n } t_{ n } - 1 }{ d } \right]  } L 
\left( \S_{ \B }^{\circ}   , t_{ 2 } - a_{ 21 } m , \dots , t_{ n+1 } - a_{ 
n+1,1 } m \right)  \end{equation}
and
   \begin{equation}\label{clos} L \left( \overline{\S_{ \A }} , \t \right) 
= \sum_{ m = \left[ \frac{ t_{ 1 } - 1 }{ a_{ 11 }  }  \right] + 1 }^{ 
\left[ \frac{ c_{ 2 } t_{ 2 } + \dots + c_{ n } t_{ n } }{ d } \right]  } L 
\left( \overline{\S_{ \B }} , t_{ 2 } - a_{ 21 } m , \dots , t_{ n+1 } - 
a_{ n+1,1 } m \right)  \ , \end{equation}
respectively, where
   \[ \B = \left( \begin{array}{rcl} a_{ 22 } & \dots & a_{ 2n } \\
                   & \vdots \\
                   a_{ n+1,2 } & \dots & a_{ n+1,n } \end{array}  \right) 
\in M_{ n \times (n-1) } ( \Z ) \ .  \]
Note that if we start with some $ \t \in \Z^{ n+1 } $ which satisfies 
Definition \ref{def}, then the
dilation parameters for $ \S_{ \B } $ in (\ref{int}) and (\ref{clos}) will 
ensure well-definedness
of the lattice point count operators.
$ L \left( \S_{ \B }^{\circ}   , \t \right) $ and $ L \left( \overline{\S_{ 
\B }} , \t \right) $ are, by induction hypothesis, quasipolynomials
satisfying the reciprocity law (\ref{rec}). Hence, by Lemma 
\ref{quasilemma}, $ L \left( \S_{ \A }^{\circ}  , \t \right) $ and $ L 
\left( \overline{\S_{ \A }} , \t \right) $
are also quasipolynomials. Note that we again use (\ref{sum}) to define 
these expressions for
all $ \t \in \Z^{ n+1 } $. Furthermore,
\begin{eqnarray*} &\mbox{}& L \left( \S_{ \A }^{\circ}  , -\t \right) = 
\sum_{ m = \left[ \frac{ -t_{ 1 } }{ a_{ 11 }  }  \right] + 1 }^{ \left[ 
\frac{ -c_{ 2 } t_{ 2 } - \dots - c_{ n } t_{ n } - 1 }{ d } \right]  } L 
\left( \S_{ \B }^{\circ}   , -t_{ 2 } - a_{ 21 } m , \dots , -t_{ n+1 } - 
a_{ n+1,1 } m \right) \\
                   &\mbox{}& \qquad \stackrel{(\ref{rec}),(\ref{sum})}{=} - 
\sum_{ \left[ \frac{ -c_{ 2 } t_{ 2 } - \dots - c_{ n } t_{ n } - 1 }{ d } 
\right] + 1 }^{ \left[ \frac{ -t_{ 1 } }{ a_{ 11 }  }  \right] } (-1)^{ n-1 
}  L \left( \overline{\S_{ \B }} , t_{ 2 } + a_{ 21 } m , \dots , t_{ n+1 } 
+ a_{ n+1,1 } m \right) \\
                   &\mbox{}& \qquad \stackrel{(\ref{gi})}{=} (-1)^{ n } 
\sum_{ m = - \left[ \frac{ c_{ 2 } t_{ 2 } + \dots + c_{ n } t_{ n } }{ d } 
\right] }^{ - \left[ \frac{ t_{ 1 } - 1 }{ a_{ 11 }  }  \right] - 1 } L 
\left( \overline{\S_{ \B }} , t_{ 2 } + a_{ 21 } m , \dots , t_{ n+1 } + 
a_{ n+1,1 } m \right)  \\
                   &\mbox{}& \qquad = (-1)^{ n } \sum_{ m = \left[ \frac{ 
t_{ 1 } - 1 }{ a_{ 11 }  }  \right] + 1 }^{ \left[ \frac{ c_{ 2 } t_{ 2 } + 
\dots + c_{ n } t_{ n } }{ d } \right] } L \left( \overline{\S_{ \B }} , 
t_{ 2 } - a_{ 21 } m , \dots , t_{ n+1 } - a_{ n+1,1 } m \right)  \\
                   &\mbox{}& \qquad = (-1)^{ n } L \left( \overline{\S_{ \A 
}} , \t \right) \ . \end{eqnarray*}
\hfill {} $\Box$


\section{\normalsize Some remarks and an example}
An obvious generalization of Theorem \ref{multithm} would be a similar 
statement for arbitrary
rational polytopes (with any number of facets). However, it is not even 
clear how to phrase
conditions on $ \t $ in the definition of a 'vector-dilated polytope', 
since the number
of facets/vertices changes for different values of $ \t $.

Another variation of the idea of vector-dilating a polytope is to dilate 
the {\it vertices}
by certain factors, instead of the facets. This would most certainly 
require completely different
methods as the ones used in this paper.

It is, finally, of interest to compute precise formulas (that is, the 
coefficients of the quasipolynomials)
for $ L \left( \S_{ \A }^{\circ}  , \t \right) $ and $ L \left( 
\overline{\S_{ \A }} , \t \right) $,
corresponding to the various existing formulas for $ L \left( {\cal 
P}^{\circ}  , t \right) $ and $ L \left( \overline{\cal P} , t \right) $.

To illustrate this, we will compute $ L \left( \overline{\S_{ \A }} , \t 
\right) $ for
a two-dimensional rectangular rational triangle, namely,
   \[ \S_{ \A } = \left\{ \x \in \R^{ 2 } : \ \begin{array}{ccccc} a_{ 1 } 
x_{ 1 } & & & \geq & 1 \\
                                                                   & & a_{ 
2 } x_{ 2 } & \geq & 1 \\
                                                                   c_{ 1 } 
x_{ 1 } & + & c_{ 2 } x_{ 2 } & \leq & 1
                                              \end{array} \right\} \ . \]
Here, $ a_{ 1 } , a_{ 2 } , c_{ 1 } , c_{ 2 } $ are positive integers;
we may also assume that $ c_{ 1 } $ and $ c_{ 2 } $ are relatively prime.
To derive a formula for $ L \left( \overline{\S_{ \A }} , \t \right) $ we 
use the methods
introduced in \cite{beck}. Similarly as in that paper, we can interpret
   \[ L \left( \overline{\S_{ \A }} , \t \right) = \# \left\{ ( m_{ 1 } , 
m_{ 2 } ) \in \Z^{ 2 } : \ \begin{array}{ccccc} a_{ 1 } m_{ 1 } & & & \geq 
& t_{ 1 }  \\
                                                                   & & a_{ 
2 } m_{ 2 } & \geq & t_{ 2 }  \\
                                                                   c_{ 1 } 
m_{ 1 } & + & c_{ 2 } m_{ 2 } & \leq & t_{ 3 }
                                              \end{array} \right\} \]
as the Taylor coefficient of $ z^{ t_{ 3 }  } $ of the function
\begin{eqnarray*} &\mbox{}& \left( \sum_{ m_{ 1 } \geq \left[ \frac{ t_{ 1 
} - 1 }{ a_{ 1 }  }  \right] + 1 } z^{ c_{ 1 } m_{ 1 }  }  \right) \left( 
\sum_{ m_{ 2 } \geq \left[ \frac{ t_{ 2 } - 1 }{ a_{ 2 }  }  \right] + 1 } 
z^{ c_{ 2 } m_{ 2 }  }  \right) \left( \sum_{ k \geq 0 } z^{ k }  \right) \\
                   &\mbox{}& \quad = \frac{ z^{ \left( \left[ \frac{ t_{ 1 
} - 1 }{ a_{ 1 }  }  \right] + 1 \right) c_{ 1 }  }  }{ 1 - z^{ c_{ 1 
}  }  } \frac{ z^{ \left( \left[ \frac{ t_{ 2 } - 1 }{ a_{ 2 }  }  \right] 
+ 1 \right) c_{ 2 }  }  }{ 1 - z^{ c_{ 2 }  }  } \frac{ 1 }{ 1 - z }   \ . 
\end{eqnarray*}
Equivalently,
   \begin{equation}\label{residue} L \left( \overline{\S_{ \A }} , \t 
\right) = \mbox{Res} \left( \frac{ z^{ e_{ 1 } + e_{ 2 } - t_{ 3 }  - 1 } 
}{ \left( 1 - z^{ c_{ 1 } }  \right) \left( 1 - z^{ c_{ 2 } }  \right) 
\left( 1 - z \right) } , z=0 \right) \ , \end{equation}
where we introduced, for ease of notation, $ e_{ j } := \left( \left[ 
\frac{ t_{ j } - 1 }{ a_{ j }  }  \right] + 1 \right) c_{ j } $
for $ j = 1, 2 $.
If the right-hand side of (\ref{residue}) counts the number of lattice 
points in $ \S_{ \A }^{ (\t) } $,
then the remaining task is computing the other residues of
   \[ f (z) := \frac{ z^{ e_{ 1 } + e_{ 2 } - t_{ 3 }  - 1 } }{ \left( 1 - 
z^{ c_{ 1 } }  \right) \left( 1 - z^{ c_{ 2 } }  \right) \left( 1 - z 
\right) } \ , \]
and use the residue theorem for the sphere $\C \cup \left\{ \infty \right\} 
$. Besides at 0, $ f $ has poles at
all $ c_{1}, c_{2} $'th roots of unity; note that if we start with a $ \t $ 
which satisfies
Definition \ref{def} then Res($ f(z), z = \infty $) = 0.

The residue at $ z=1 $ can be easily calculated as
   \begin{eqnarray*} &\mbox{}& \mbox{Res} \Bigl( f (z) , z=1 \Bigr) = 
\mbox{Res} \Bigl( e^{ z }  f (e^{ z } ) , z=0 \Bigr)  \\
                     &\mbox{}& \quad = - \frac{ 1 }{ 2c_{ 1 } c_{ 2 }  } 
\left( e_{ 1 } + e_{ 2 } - t_{ 3 } \right)^{ 2 } + \frac{ 1 }{ 2 } \left( 
e_{ 1 } + e_{ 2 } - t_{ 3 } \right) \left( \frac{ 1 }{ c_{ 1 }  } + \frac{ 
1 }{ c_{ 2 }  } + \frac{ 1 }{ c_{ 1 } c_{ 2 }  }  \right) \\
                     &\mbox{}& \quad \qquad - \frac{ 1 }{ 4 } \left( 1 + 
\frac{ 1 }{ c_{ 1 }  } + \frac{ 1 }{ c_{ 2 }  }  \right) - \frac{ 1 }{ 12 } 
\left( \frac{ c_{ 1 }  }{ c_{ 2 }  } + \frac{ c_{ 2 }  }{ c_{ 1 }  } + 
\frac{ 1 }{ c_{ 1 } c_{ 2 }  }  \right) \ . \end{eqnarray*}
It remains to compute the residues at the nontrivial roots of unity. Let $ 
\lambda^{ c_{ 1 } } = 1 \not= \lambda $.
Then
   \begin{eqnarray*} &\mbox{}& \mbox{Res} \Bigl( f (z) , z=\lambda \Bigr) = 
\frac{ \lambda^{ e_{ 2 } - t_{ 3 } - 1}  }{ \left( 1 - \lambda^{ c_{ 2 
}  }  \right) \left( 1 - \lambda \right)  } \ \mbox{Res} \left( \frac{ 1 }{ 
1 - \lambda^{ c_{ 1 }  }  }  , z=\lambda \right) \\
                     &\mbox{}& \quad = - \frac{ \lambda^{ e_{ 2 } - t_{ 3 } 
}  }{ c_{ 1 } \left( 1 - \lambda^{ c_{ 2 }  }  \right) \left( 1 - \lambda 
\right)  } \ . \end{eqnarray*}
Adding up all the nontrivial $ c_{ 1 } $'th roots of unity, we obtain
   \[ \sum_{ \lambda^{ c_{ 1 } } = 1 \not= \lambda  }  \mbox{Res} \Bigl( f 
(z) , z=\lambda \Bigr) = - \frac{ 1 }{ c_{ 1 }  }  \sum_{ \lambda^{ c_{ 1 } 
} = 1 \not= \lambda  } \frac{ \lambda^{ e_{ 2 } - t_{ 3 } }  }{ \left( 1 - 
\lambda^{ c_{ 2 }  }  \right) \left( 1 - \lambda \right)  } \ , \]
a special case of a {\it Fourier-Dedekind sum}, which already occurred in 
\cite{bdr}.
In fact, in the same paper we derived, by means of finite Fourier series,
   \[ \frac{ 1 }{ c_{ 1 }  }  \sum_{ \lambda^{ c_{ 1 } } = 1 \not= 
\lambda  } \frac{ \lambda^{ t }  }{ \left( 1 - \lambda^{ c_{ 2 
}  }  \right) \left( 1 - \lambda \right)  } = \sum_{ k=0 }^{ c_{ 1 } - 1 } 
\left( \left( \frac{ - c_{ 2 } k - t }{ c_{ 1 }  }  \right)  \right) \left( 
\left( \frac{ k }{ c_{ 1 }  }  \right)  \right) - \frac{ 1 }{ 4 c_{ 1 }  } 
\ , \]
where $ (( x )) = x - [x] - 1/2 $ is a sawtooth function (differing 
slightly from
the one appearing in the classical Dedekind sums). The expression on the 
right is, up to a
trivial term, a special case of a {\it Dedekind-Rademacher sum} 
(\cite{dieter}, \cite{meyer},
\cite{rademacher}). Hence,
   \[ \sum_{ \lambda^{ c_{ 1 } } = 1 \not= \lambda  } \mbox{Res} \Bigl( f 
(z) , z=\lambda \Bigr) = - \sum_{ k=0 }^{ c_{ 1 } - 1 } \left( \left( 
\frac{ t_{ 3 } - e_{ 2 } - c_{ 2 } k }{ c_{ 1 }  }  \right)  \right) \left( 
\left( \frac{ k }{ c_{ 1 }  }  \right)  \right) + \frac{ 1 }{ 4 c_{ 1 }  } 
\ , \]
and, similarly, for the nontrivial $ c_{ 2 } $'th roots of unity
   \[ \sum_{ \mu^{ c_{ 2 } } = 1 \not= \mu } \mbox{Res} \Bigl( f (z) , 
z=\mu \Bigr) = - \sum_{ k=0 }^{ c_{ 2 } - 1 } \left( \left( \frac{ t_{ 3 } 
- e_{ 1 } - c_{ 1 } k }{ c_{ 2 }  }  \right)  \right) \left( \left( \frac{ 
k }{ c_{ 2 }  }  \right)  \right) + \frac{ 1 }{ 4 c_{ 2 }  } \ . \]
The residue theorem allows us now to rewrite (\ref{residue}) as
   \begin{eqnarray*} &\mbox{}& L \left( \overline{\S_{ \A }} , \t \right) = 
\frac{ 1 }{ 2 c_{ 1 } c_{ 2 }  } \left( e_{ 1 } + e_{ 2 } - t_{ 3 } 
\right)^{ 2 } - \frac{ 1 }{ 2 } \left( e_{ 1 } + e_{ 2 } - t_{ 3 } \right) 
\left( \frac{ 1 }{ c_{ 1 }  } + \frac{ 1 }{ c_{ 2 }  } + \frac{ 1 }{ c_{ 1 
} c_{ 2 }  }  \right) \\
                     &\mbox{}& \quad \qquad + \frac{ 1 }{ 4 } + \frac{ 1 }{ 
12 } \left( \frac{ c_{ 1 }  }{ c_{ 2 }  } + \frac{ c_{ 2 }  }{ c_{ 1 }  } + 
\frac{ 1 }{ c_{ 1 } c_{ 2 }  }  \right) + \sum_{ k=0 }^{ c_{ 1 } - 1 } 
\left( \left( \frac{ t_{ 3 } - e_{ 2 } - c_{ 2 } k }{ c_{ 1 
}  }  \right)  \right) \left( \left( \frac{ k }{ c_{ 1 
}  }  \right)  \right) \\
                     &\mbox{}& \quad \qquad + \sum_{ k=0 }^{ c_{ 2 } - 1 } 
\left( \left( \frac{ t_{ 3 } - e_{ 1 } - c_{ 1 } k }{ c_{ 2 
}  }  \right)  \right) \left( \left( \frac{ k }{ c_{ 2 
}  }  \right)  \right) \ . \end{eqnarray*}
To see the quasipolynomial character better, we substitute back the 
expressions for $ e_{ 1 } $
and $ e_{ 2 } $, and write $ [x] = x - ((x)) - 1/2 $ for the greatest 
integer function.
After a somewhat tedious calculation, we obtain
   \begin{eqnarray*} &\mbox{}& L \left( \overline{\S_{ \A }} , \t \right) = 
\frac{ c_{ 1 }  }{ 2 a_{ 1 }^{ 2 } c_{ 2 }   } t_{ 1 }^{ 2 } + \frac{ c_{ 2 
}  }{ 2 a_{ 2 }^{ 2 } c_{ 1 }   } t_{ 2 }^{ 2 } + \frac{ 1 }{ 2 c_{ 1 } c_{ 
2 }  } t_{ 3 }^{ 2 } + \frac{ 1 }{ a_{ 1 } a_{ 2 }  } t_{ 1 } t_{ 2 } - 
\frac{ 1 }{ a_{ 1 } c_{ 2 }  } t_{ 1 } t_{ 3 } - \frac{ 1 }{ a_{ 2 } c_{ 1 
}  } t_{ 2 } t_{ 3 } \\
                     &\mbox{}& \qquad \qquad \qquad + \nu_{ 1 } ( \t ) \ 
t_{ 1 } + \nu_{ 2 } ( \t ) \ t_{ 2 } + \nu_{ 3 } ( \t ) \ t_{ 3 } + \nu_{ 0 
} ( \t ) \ ,  \end{eqnarray*}
where
   \begin{eqnarray*} &\mbox{}& \nu_{ 1 } ( \t ) = - \frac{ c_{ 1 }  }{ a_{ 
1 }^{ 2 } c_{ 2 }   } \left( 1 + \left( \left( \frac{ t_{ 1 } - 1 }{ a_{ 1 
}  }  \right)  \right)  \right) - \frac{ 1 }{ a_{ 1 }  } \left( \left( 
\frac{ t_{ 2 } - 1 }{ a_{ 2 }  }  \right)  \right) - \frac{ 1 }{ a_{ 1 } 
a_{ 2 }  } - \frac{ 1 }{ 2 a_{ 1 } c_{ 2 }  } \\
                     &\mbox{}& \nu_{ 2 } ( \t ) = - \frac{ c_{ 2 }  }{ a_{ 
2 }^{ 2 } c_{ 1 }   } \left( 1 + \left( \left( \frac{ t_{ 2 } - 1 }{ a_{ 2 
}  }  \right)  \right)  \right) - \frac{ 1 }{ a_{ 2 }  } \left( \left( 
\frac{ t_{ 1 } - 1 }{ a_{ 1 }  }  \right)  \right) - \frac{ 1 }{ a_{ 1 } 
a_{ 2 }  } - \frac{ 1 }{ 2 a_{ 2 } c_{ 1 }  } \\
                     &\mbox{}& \nu_{ 3 } ( \t ) = \frac{ 1 }{ a_{ 1 } c_{ 2 
}  } + \frac{ 1 }{ a_{ 2 } c_{ 1 }  } + \frac{ 1 }{ 2 c_{ 1 } c_{ 2 }  } + 
\frac{ 1 }{ c_{ 2 }  } \left( \left( \frac{ t_{ 1 } - 1 }{ a_{ 1 
}  }  \right)  \right) + \frac{ 1 }{ c_{ 1 }  } \left( \left( \frac{ t_{ 2 
} - 1 }{ a_{ 2 }  }  \right)  \right) \\
                     &\mbox{}& \nu_{ 0 } ( \t ) = - \frac{ 1 }{ 4 c_{ 1 
}  } - \frac{ 1 }{ 4 c_{ 2 }  } + \frac{ 1 }{ a_{ 1 } a_{ 2 }  } + \frac{ 1 
}{ 2 a_{ 1 } c_{ 2 }  } + \frac{ 1 }{ 2 a_{ 2 } c_{ 1 }  } + \frac{ 1 }{ 12 
c_{ 1 } c_{ 2 }  } - \frac{ c_{ 1 }  }{ 24 c_{ 2 }  } - \frac{ c_{ 2 }  }{ 
24 c_{ 1 }  } \\
                     &\mbox{}& \quad \qquad + \frac{ c_{ 1 }  }{ 2 a_{ 1 
}^{ 2 } c_{ 2 }   } + \frac{ c_{ 2 }  }{ 2 a_{ 2 }^{ 2 } c_{ 1 }   } + 
\left( \left( \frac{ t_{ 1 } - 1 }{ a_{ 1 }  }  \right)  \right) \left( 
\frac{ 1 }{ a_{ 2 }  } + \frac{ 1 }{ 2 c_{ 2 }  } + \frac{ c_{ 1 }  }{ a_{ 
1 } c_{ 2 }  }  \right) \\
                     &\mbox{}& \quad \qquad + \left( \left( \frac{ t_{ 2 } 
- 1 }{ a_{ 2 }  }  \right)  \right) \left( \frac{ 1 }{ a_{ 1 }  } + \frac{ 
1 }{ 2 c_{ 1 }  } + \frac{ c_{ 2 }  }{ a_{ 2 } c_{ 1 }  }  \right) + \frac{ 
c_{ 1 }  }{ 2 c_{ 2 }  } \left( \left( \frac{ t_{ 1 } - 1 }{ a_{ 1 
}  }  \right)  \right)^{ 2 } \\
                     &\mbox{}& \quad \qquad + \frac{ c_{ 2 }  }{ 2 c_{ 1 
}  } \left( \left( \frac{ t_{ 2 } - 1 }{ a_{ 2 }  }  \right)  \right)^{ 2 } 
+ \left( \left( \frac{ t_{ 1 } - 1 }{ a_{ 1 }  }  \right)  \right) \left( 
\left( \frac{ t_{ 2 } - 1 }{ a_{ 2 }  }  \right)  \right) \\
                     &\mbox{}& \quad \qquad + \sum_{ k=0 }^{ c_{ 1 } - 1 } 
\left( \left( \frac{ t_{ 3 } }{ c_{ 1 }  } - \frac{ t_{ 2 } - 1 }{ a_{ 2 } 
c_{ 1 }  } + \frac{ 1 }{ c_{ 1 }  } \left( \left( \frac{ t_{ 2 } - 1 }{ a_{ 
2 }  }  \right)  \right) - \frac{ 1 }{ 2 c_{ 1 }  } - \frac{ c_{ 2 } k }{ 
c_{ 1 }  }  \right)  \right) \left( \left( \frac{ k }{ c_{ 1 
}  }  \right)  \right) \\
                     &\mbox{}& \quad \qquad + \sum_{ k=0 }^{ c_{ 2 } - 1 } 
\left( \left( \frac{ t_{ 3 } }{ c_{ 2 }  } - \frac{ t_{ 1 } - 1 }{ a_{ 1 } 
c_{ 2 }  } + \frac{ 1 }{ c_{ 2 }  } \left( \left( \frac{ t_{ 1 } - 1 }{ a_{ 
1 }  }  \right)  \right) - \frac{ 1 }{ 2 c_{ 2 }  } - \frac{ c_{ 1 } k }{ 
c_{ 2 }  }  \right)  \right) \left( \left( \frac{ k }{ c_{ 2 
}  }  \right)  \right) \ . \end{eqnarray*}
As a final remark, we note that this formula enables us to compute the 
number of
lattice points inside {\it any} rational polygon: Any two-dimensional 
polytope can be
written as a virtual decomposition of rectangles (which are easy to deal 
with) and the right-angled triangles
discussed above. Moreover, if the polygon has rational vertices, so do all 
these 'pieces'.


\vspace{5mm}

{\bf Acknowledgements}. I am grateful to Boris Datskovsky, Sinai Robins, 
and Bob Styer
for corrections and helpful suggestions on earlier versions of this paper, 
and to
Tendai Chitewere for invaluable moral support.

\small
\nocite{*}
\addcontentsline{toc}{subsubsection}{References}
\bibliography{thesis}
\bibliographystyle{alpha}

\end{document}